\documentclass[a4paper,12pt]{amsart}
\usepackage{amssymb}
\usepackage{amsmath}
\usepackage{ifthen}
\nonstopmode \numberwithin{equation}{section}
\newtheorem{thm}[equation]{Theorem}
\newtheorem{cor}[equation]{Corollary}
\newtheorem{lem}[equation]{Lemma}
\newtheorem{prop}[equation]{Proposition}

\newenvironment{pf}[1][]{%
 \vskip 3mm
 \noindent
 \ifthenelse{\equal{#1}{}}%
  {{\slshape Proof. }}%
  {{\slshape #1.} }%
 }%
{\qed\bigskip}

\begin{document}
\bibliographystyle{amsplain}
\newcommand{\A}{{\mathcal A}}
\newcommand{\B}{{\mathcal B}}
\newcommand{\UCV}{{\mathcal UCV}}
\newcommand{\UST}{{\mathcal UST}}
\newcommand{\RR}{{\mathcal R}}
\newcommand{\es}{{\mathcal S}}
\newcommand{\IR}{{\mathbb R}}
\newcommand{\IC}{{\mathbb C}}
\newcommand{\IN}{{\mathbb N}}
\newcommand{\K}{{\mathcal K}}
\newcommand{\uhp}{{\mathbb H}}
\newcommand{\Z}{{\mathbb Z}}
\newcommand{\N}{{\mathcal N}}
\newcommand{\M}{{\mathcal M}}
\newcommand{\SCC}{{\mathcal{SCC}}}
\newcommand{\CC}{{\mathcal C}}
\newcommand{\st}{{\mathcal{SS}}}
\newcommand{\D}{{\mathbb D}}
\newcommand{\remark}{\vskip .3cm \noindent {\sl Remark.} \@}
\newcommand{\remarks}{\vskip .3cm \noindent {\sl Remarks.} \@}
\newcommand{\ucv}{{\operatorname{UCV}}}
\def\be{\begin{equation}}
\def\ee{\end{equation}}
\newcommand{\bee}{\begin{enumerate}}
\newcommand{\eee}{\end{enumerate}}
\newcommand{\pays}{\!\!\!\!}
\newcommand{\pay}{\!\!\!}
\newcommand{\blem}{\begin{lem}}
\newcommand{\elem}{\end{lem}}
\newcommand{\bthm}{\begin{thm}}
\newcommand{\ethm}{\end{thm}}
\newcommand{\bcor}{\begin{cor}}
\newcommand{\ecor}{\end{cor}}
\newcommand{\beg}{\begin{example}}
\newcommand{\eeg}{\end{example}}
\newcommand{\begs}{\begin{examples}}
\newcommand{\eegs}{\end{examples}}
\newcommand{\bdefe}{\begin{defn}}
\newcommand{\edefe}{\end{defn}}
\newcommand{\bprob}{\begin{prob}}
\newcommand{\eprob}{\end{prob}}
\newcommand{\bcon}{\begin{conj}}
\newcommand{\econ}{\end{conj}}
\newcommand{\bprop}{\begin{prop}}
\newcommand{\eprop}{\end{prop}}
\newcommand{\bpf}{\begin{pf}}
\newcommand{\epf}{\end{pf}}
\newcommand{\ba}{\begin{array}}
\newcommand{\ea}{\end{array}}
\newcommand{\beq}{\begin{eqnarray}}
\newcommand{\beqq}{\begin{eqnarray*}}
\newcommand{\eeq}{\end{eqnarray}}
\newcommand{\eeqq}{\end{eqnarray*}}

\title[New Hohlov Type Integral Operator involving Clausen's  Hypergeometric Functions]
{New Hohlov Type Integral Operator involving Clausen's  Hypergeometric Functions}

\author[K.Chandrasekran]{K. Chandrasekran}
\address{K. Chandrasekran \\ Department of Mathematics \\ MIT Campus, Anna University \\ Chennai 600 044, India}
\email{kchandru2014@gmail.com}

\author[D. J. Prabhakaran ]{D.J. Prabhakaran}
\address{D.J. Prabhakaran \\ Department of Mathematics \\ MIT Campus, Anna University \\ Chennai 600 044, India}
\email{asirprabha@gmail.com}

\subjclass[2000]{30C45}
\keywords{Clausen's  Hypergeometric Function, Univalent Functions, Starlike Functions, Convex Functions, Uniformly Starlike Functions and Uniformly Convex Functions.}

\begin{abstract}
We consider the integral operator $\mathcal{I}^{a,b,c}_{d,e}(f)(z)$ involving Clausen's Hypergeometric Function by means of convolution introduced by Chandrasekran and Prabhakaran for investigation. The conditions on the parameters $ a,b, c$ are determined using the integral operator $\mathcal{I}^{a,\frac{b}{2},\frac{b+1}{2}}_{\frac{c}{2}, \frac{c+1}{2}}(f)(z)$  to study the geometric properties of Clausen's Hypergeometric Function for various subclasses of univalent functions.
\end{abstract}
\maketitle

\section{Introduction and preliminaries}

Let $\A$ denote the family of analytic functions $f$ of the form
\beq\label{inteq0}
f(z)= z+\sum_{n=2}^{\infty}\, a_n\,z^n
\eeq with $f(0)=0$ and $f^{\prime}(0)=1$ in the open unit disc $\D =\{z:\, |z|<1\}$ of the complex plane.\\

Let us denote by $\es,$ the class of all normalised functions that are analytic and univalent in $\D;$ that is $ {\es}=\{f(z) \in \A: f(z)  \text{ is univalent in\,} \ \D\}.$ If $\displaystyle f(z)$ defined by (\ref{inteq0}) is belongs to $\es$, then,
\beq\label{inteq00}
|a_n|\leq n,\text{ for}\, n\geq 2.
\eeq

Various subclasses of $\es$ are characterized by their geometric properties. Some important subclasses of $\es$ are the class $\CC$ of all convex functions, the class $\es^{*}$ of all functions starlike with respect to the origin and the class $\K$ of all close-to-convex functions.(For more details refer \cite{Peter-L-Duren-book-1983,A-W-Goodman-1983-book}.\\

Robertson \cite{Robertson-1936} introduced the class $\es^{*}{(\alpha)},$ of all starlike functions of order $\alpha $ and the class $\CC{(\alpha)}$ of all convex functions of order $\alpha.$ A function in  $f(z) \in \es$ is a starlike function of order $\alpha$ in $\D$ if and only if  ${ \Re\,} \displaystyle \left( \frac{zf^{\prime}(z)}{f(z)} \right) > \alpha, \text{ for all } z\in \D$ where $0\leq \alpha <1.$ A function in  $f(z) \in \es$ is a convex function of order $\alpha$ in $\D$ if and only if $ \displaystyle { \Re\, } \left(1+ \frac{zf^{\prime\prime}(z)}{f^{\prime}(z)} \right) >  \alpha \text{ for all } z \in \D,$ where $0\leq \alpha <1.$\\
%

%
Our main focus is on the following subclasses of starlike and convex functions that are defined as follows. For $\lambda >0,\, \es^{*}_{\lambda},\,  $ is  defined as follows $$\es^{*}_{\lambda}\, =\, \left\{f(z)\in \A \, |\, \left|\frac{zf'(z)}{f(z)}-1\right| \, < \, \lambda,\, z\in \D \right\}$$ and a sufficient condition for which the function $f(z)$ to be in $\es^{*}_{\lambda}$ is
\beq\label{inteq2}
\displaystyle \sum_{n=2}^{\infty}(n+\lambda-1)|a_n| \leq \lambda.
\eeq

For $\lambda >0,$ $\CC_{\lambda}$ is defined as follows $\CC_{\lambda}=\left\{f(z)\in \A\, |\, zf'(z)\in \es^{*}_{\lambda}\right\}$ and a sufficient condition for which the function $f(z)$ to be in $\CC_{\lambda}$, is as follows:
\beq\label{inteq}
\displaystyle \sum_{n=2}^{\infty}\, n\, (n+\lambda-1)|a_n| \leq \lambda.
\eeq

Goodman \cite{Good-1991-Ann-PM, Good-1991-JMAA} introduced the two new class, the class of all uniformly convex functions denoted by ${\UCV}$ and the class of all uniformly starlike functions denoted by ${\UST}$. He gave an analytic criterion for a function to be uniformly convex and u   niformly starlike. 
Let $f(z)$ be a function of the form (\ref{inteq0}). Then $f$ is in ${\UCV}$ if and only if $\displaystyle\Re \left(1 + \frac{(z - \zeta)f^{\prime\prime}(z)} {f^{\prime}(z)}\right) > 0,$ for all $z, \ \zeta \in \D.$ Let $f(z)$ be a function of the form (\ref{inteq0}). Then $f$ is in ${\UST}$ if and only if $\displaystyle\Re \left( \frac{(z - \zeta)f^{\prime}(z)} {f(z)-f(\zeta)} \right) > 0,$ for all $z, \  \zeta \in \D.$ Note that, for $\zeta = 0,$ the class ${\UCV}$ coincides with the class $\CC$ and the class ${\UST}$ coincides with the class ${\es}^{*}.$
Subsequently, R{\o}nning \cite{Ronn-1993-Proc-ams} and Ma and Minda \cite{Ma-Minda-1992-Ann-PM} independently gave the one variable analytic characterization of the class $\UCV$.\\

A sufficient condition for the function of form (\ref{inteq0}) to belong to $\UCV$ \cite{Subram-Murugu-1995}, is given by
\beq\label{lem4eq1}
\sum_{n=2}^{\infty}\,n\, (2n-1)|a_n|\leq 1.
\eeq

The  subclass $\es_p$ of starlike functions introduced by R{\o}nning \cite{Ronn-1993-Proc-ams} is defined as
\beqq
\displaystyle \es_p = \{ F \in \es^{*}|F(z)=zf'(z),\, f(z) \in \UCV\}.
\eeqq

A sufficient condition for a function $f(z)$ of form (\ref{inteq0}) to belong to $\es_p$ is given by
 \beq \label{lem2eq1}
\sum_{n=2}^{\infty}(2n-1)|a_n|\leq 1.
\eeq

The above result was proved for more general case $\es_p(\alpha)$ in \cite{Subra-Sudharsan-1998}.\\

Let for $\beta <1$,
\beqq
  {\RR}(\beta)=\{f(z)\in {\A}:  \exists \  \phi \in {\IR} / { Re}\, \left( e^{i\phi}(f'(z)-\beta)\right)>0, \quad z\in\D \}.
\eeqq

Note that when $\beta \ge 0$ we have ${\RR}(\beta) \subset {\es}$ and for $\beta <0 ,\ \  {\RR}(\beta) $ contains also non univalent functions. This class has been widely used to study certain integral transforms. Further \cite{Anbu-Parva-2000, Parva-prabha-2001-Far-East} and the reference therein. Suppose that $\displaystyle f(z)$ defined by (\ref{inteq0}) is in the class $ \RR(\beta)$. Then, by \cite{MacGregor-1962-Trans-ams}, we have
\beq\label{inteq3}
|a_n|\leq \frac{2(1-\beta)}{n},\, n \geq 2.
\eeq

The convolution of two functions $\displaystyle f(z)= z+\sum_{n=2}^{\infty}\, a_n\,z^n $ and  $ \displaystyle g(z)= z+\sum_{n=2}^{\infty}\, b_n\,z^n $ (both $f(z)$ and $g(z)$ are analytic functions in $\D$ ) is defined by $\displaystyle f(z)*g(z)= z+\sum_{n=2}^{\infty}\, a_n\,b_n\, z^n.$\\

For any non-zero complex variable $a$, the ascending factorial notation (or Pochhammer symbol) is defined as
$(a)_0\,=\,1,\, {\rm and }\, (a)_n\,=\,a(a+1)\cdots (a+n-1),\, {\rm for}\, n\,=\,1,2,3,\cdots.$\\

In literature, there are many articles connecting geometric function theory and Gauss hypergeometric function. But there are very few papers on Clausen’s hypergeometric function. This function is of greater interest, since the famous Bieberbach conjecture has prove by Louis De Branges in 1984 using the Askey - Gasper inequalities for Jacobi Polynomials which involves Clausen’s series.
$$ _3F_2\left(\begin{array}{c}
                       k+\frac{1}{2},\, k-n,\, n+k+2 \\
                       2k+1, \, k+\frac{3}{2}
                     \end{array}; e^{-t}\right).$$

The Clausen's series can be obtained by squaring the gaussian hypergeometric function. i.e.,
 $\displaystyle [_2F_1(a,b;a+b+1/2;z)]^2 = \, _3F_2(2a,2b,a+b; 2a+2b,a+b+1/2;z).$\\
%

The Clausen's hypergeometric function $_3F_2(a,b,c;d,e;z)$ is defined as
\beqq
_3F_2(a,b,c;d,e;z)=\sum_{n=0}^{\infty}\frac{(a)_n(b)_n(c)_n}{(d)_n(e)_n(1)_n}z^n, \quad |z| < 1
\eeqq
where $ a,b,c,d,e\in \IC$  with  $d,\, e\, \neq 0,-1,-2,-3\cdots$.\\

The Clausen's hypergeometric function $_3F_2(a,b,c;d,e;z) (:=\,  _4F_3(a,b,c,a_4;d,e,a_4;z))$ have been studied by only few authors. In particular \cite{Ponnu-saba-1997}, Ponnusamy and Sabapathy considered the generalized hypergeometric functions and tried to find condition on the parameter so that $z\,_3F_2(a,b,c;d,e;z)$ has some geometric properties.\\

In \cite{Chandru-prabha-2019}, Chandrasekran and Prabhakaran have introduced an integral operator and derived the geometric properties for the Clausen's Hypergeometric series $z\, _3F_2(a,b,c;b+1,c+1;z)$, in which, the numerator and denominator parameters differs by arbitrary negative integers. Also to determine the conditions on the parameters $a,\,b,\,c,$ such that, the hypergeometric function $z\,_3F_2(a,b,c;1+a-b,\, 1+a-c)$ associated with the Dixon's summation formula or its equivalence has the admissibility to be the classes like $\mathcal{S}^{*}_{\lambda}, \, \mathcal{C}_{\lambda}$,\, $\mathcal{UCV}$, and $\es_p$ is still it is an open problem.\\

In 2006, Driver and Johnston  \cite{Driver-Johnston-2006} derived a summation formula for Clausen's hypergeometric function in terms of Gaussian hypergeometric function. We recall their summation formula as follows:
\begin{sloppypar}
\beq\label{inteq6}
_3F_2\left(a,\frac{b}{2},\frac{b+1}{2};\frac{c}{2}, \frac{c+1}{2}; 1\right)&=& \frac{\Gamma(c)\,\Gamma(c-a-b)}{\Gamma(c-a)\,\Gamma(c-b)}\, _{2}F_1(a,b;c-a;-1)
\eeq provide $Re(c)\, >\, Re(b)\, > 0$ and $Re(c-a-b)\, >\, 0.$\\
\end{sloppypar}

In this paper, we consider an integral operator $\mathcal{I}^{a,b,c}_{d,e}(f)(z)$ introduced by Chandrasekran and Prabhakaran \cite{Chandru-prabha-2019} and to study the various geometric properties of the above Clausen's series,  the integral operator takes the form as follows
\beq\label{inteq7}
\mathcal{I}^{a,\frac{b}{2},\frac{b+1}{2}}_{\frac{c}{2}, \frac{c+1}{2}}(f)(z) &=& z\, _3F_2\left(a,\frac{b}{2},\frac{b+1}{2};\frac{c}{2}, \frac{c+1}{2};z\right)*f(z)\\
&=& z+\sum_{n=2}^{\infty} A_n\, z^n, \, f\in \mathcal{A},\nonumber
\eeq
with $A_1=1$ and for $n > 1,$
\beq\label{inteq007}
A_n&=&\frac{(a)_{n-1}\left(\frac{b}{2}\right)_{n-1}\left(\frac{b+1}{2}\right)_{n-1}}{\left(\frac{c}{2}\right)_{n-1}\left(\frac{c+1}{2}\right)_{n-1}(1)_{n-1}}\, a_n.
\eeq
The following  Lemma is useful to prove our main results.
\begin{sloppypar}
\blem \label{ch3lem1eq1}
Let $a,\, b,\, c > 0$. Then the following is derived:
\begin{enumerate}
\item For $c > a+b+1$,\
\begin{flushleft}
$\displaystyle\sum_{n=0}^{\infty} \frac{(n+1)(a)_n\, \left(\frac{b}{2}\right)_n\, \left(\frac{b+1}{2}\right)_n }
{\left(\frac{c}{2}\right)_n\, \left(\frac{c+1}{2}\right)_n\, (1)_n}$
\end{flushleft}
\begin{eqnarray*}
 &=&  \frac{\Gamma(c)\, \Gamma(c-a-b)}{\Gamma(c-a)\,\Gamma(c-b)} \,\bigg[\left(\frac{(a)\,(b)_2\,}{(c-a)\,(c-a-b-1)\,}\right)\,\cr
 && \qquad \times  _{2}F_1(a+1,b+2;c-a+1;-1)\, +\, _{2}F_1(a,b;c-a;-1)\bigg].
\end{eqnarray*}
\item For $c >  a+b+2$,\\
\begin{flushleft}
$\displaystyle\sum_{n=0}^{\infty} \frac{(n+1)^2(a)_n\, \left(\frac{b}{2}\right)_n\, \left(\frac{b+1}{2}\right)_n }
{\left(\frac{c}{2}\right)_n\, \left(\frac{c+1}{2}\right)_n\, (1)_n}$
\end{flushleft}
\begin{eqnarray*}
 &=&  \frac{\Gamma(c)\, \Gamma(c-a-b)}{\Gamma(c-a)\,\Gamma(c-b)} \,\bigg[\left(\frac{(a)_2\,(b)_4}{(c-a)_2\,(c-a-b-2)_{2}\,}\right)\, \,\cr 
 && \qquad \times _{2}F_1(a+2,b+4;c-a+2;-1)\,\cr 
 && +\, 3\left(\frac{a\,(b)_2}{(c-a)(c-a-b-1)}\right)\, _{2}F_1(a+1,b+2;c-a+1;-1)\,\cr 
 && \qquad \,+\, _{2}F_1(a,b;c-a;-1)\bigg].
\end{eqnarray*}
\item For $ c > a+b+3$,\\
\begin{flushleft}
$\displaystyle\sum_{n=0}^{\infty} \frac{(n+1)^3(a)_n\, \left(\frac{b}{2}\right)_n\, \left(\frac{b+1}{2}\right)_n }
{\left(\frac{c}{2}\right)_n\, \left(\frac{c+1}{2}\right)_n\, (1)_n}$
\end{flushleft}
\begin{eqnarray*}
 &=&  \frac{\Gamma(c)\, \Gamma(c-a-b)}{\Gamma(c-a)\,\Gamma(c-b)} \,\bigg[\left(\frac{(a)_3\,(b)_6}{(c-a)_3\,(c-a-b-3)_{3}}\right)\,\cr
 && \qquad \times\,_{2}F_1(a+3,b+6;c-a+3;-1)\cr 
 &&+ 6 \left( \frac{(a)_2\,(b)_4}{(c-a)_2\,(c-a-b-2)_{2}}\right) _{2}F_1(a+2,b+4;c-a+2;-1)\,\cr 
 && +  7 \left(\frac{a\,(b)_2}{(c-a)\,(c-a-b-1)}\right)\,\cr 
 && \qquad \times\, _{2}F_1(a+1,b+2;c-a+1;-1)\, +\, _{2}F_1(a,b;c-a;-1)\bigg].
\end{eqnarray*}
\item For $a\neq 1,\, b\neq 1,\, 2$ and $ c > \max\{a+1, a+b-1\}$, \\
\begin{flushleft}
$\displaystyle\sum_{n=0}^{\infty} \frac{(a)_n\, \left(\frac{b}{2}\right)_n\, \left(\frac{b+1}{2}\right)_n }
{\, \left(\frac{c}{2}\right)_n\, \left(\frac{c+1}{2}\right)_n\, (1)_{n+1}}$
\end{flushleft}
\begin{eqnarray*}
&=& \bigg[\left(\frac{(c-a-1)\,(c-a-b)\,\Gamma(c)\Gamma(c-a-b)}{(a-1)\, (b-1)\, (b-2)\,\Gamma(c-b)\, \Gamma(c-a)}\right)\, \cr
&&\times \, _{2}F_1(a-1,b-2;c-a-1;-1)-\left(\frac{(c-2)_2}{(a-1)\, (b-2)_2}\right)\bigg].
\end{eqnarray*}
\end{enumerate}
\elem
\bpf (1) Using Pochhammer symbol, we can formulate
\begin{flushleft}
$\displaystyle \sum_{n=0}^{\infty} \frac{(n+1)(a)_n\, \left(\frac{b}{2}\right)_n\, \left(\frac{b+1}{2}\right)_n }
{\left(\frac{c}{2}\right)_n\, \left(\frac{c+1}{2}\right)_n\, (1)_n}$
\end{flushleft}
\begin{eqnarray*}
&=& \sum_{n=0}^{\infty} \frac{(a)_{n+1}\, \left(\frac{b}{2}\right)_{n+1}\, \left(\frac{b+1}{2}\right)_{n+1} }
{\left(\frac{c}{2}\right)_{n+1}\, \left(\frac{c+1}{2}\right)_{n+1}\, (1)_{n}} +\sum_{n=0}^{\infty} \frac{(a)_n\, \left(\frac{b}{2}\right)_n\, \left(\frac{b+1}{2}\right)_n }
{\left(\frac{c}{2}\right)_n\, \left(\frac{c+1}{2}\right)_n\, (1)_n}
\end{eqnarray*}
Using the formula (\ref{inteq6}) and using the fact that $\Gamma(a+1)= a\Gamma(a)$, the aforementioned equation reduces to\\
\begin{flushleft}
$\displaystyle \sum_{n=0}^{\infty} \frac{(n+1)(a)_n\, \left(\frac{b}{2}\right)_n\, \left(\frac{b+1}{2}\right)_n }
{\left(\frac{c}{2}\right)_n\, \left(\frac{c+1}{2}\right)_n\, (1)_n}$
\end{flushleft}
\begin{eqnarray*}
 &=& \frac{\Gamma(c)\, \Gamma(c-a-b)}{\Gamma(c-b)\,\Gamma(c-a)}\\ 
 && \times\bigg[\left(\frac{(a)\, (b)_2}{(c-a)\, (c-a-b-1)}\right)\, _{2}F_1(a+1,b+2;c-a+1;-1)\,\cr 
 && \qquad \qquad \qquad  \, +\, _{2}F_1(a,b;c-a;-1)\bigg].
\end{eqnarray*}
Hence, (1) is proved.\\

(2) Using $(n+1)^2 = n(n-1)+3n+1$, we can easily obtain that
\begin{flushleft}
$\displaystyle \sum_{n=0}^{\infty} \frac{(n+1)^2 \,(a)_n\, \left(\frac{b}{2}\right)_n\, \left(\frac{b+1}{2}\right)_n }
{\left(\frac{c}{2}\right)_n\, \left(\frac{c+1}{2}\right)_n\, (1)_n}$
\end{flushleft}
\begin{eqnarray*}
&=& \sum_{n=2}^{\infty} \frac{(a)_n\, \left(\frac{b}{2}\right)_n\, \left(\frac{b+1}{2}\right)_n }
{\left(\frac{c}{2}\right)_n\, \left(\frac{c+1}{2}\right)_n\, (1)_{n-2}}+3\,\sum_{n=1}^{\infty} \frac{(a)_n\, \left(\frac{b}{2}\right)_n\, \left(\frac{b+1}{2}\right)_n }
{\left(\frac{c}{2}\right)_n\, \left(\frac{c+1}{2}\right)_n\, (1)_{n-1}} +\sum_{n=0}^{\infty} \frac{(a)_n\, \left(\frac{b}{2}\right)_n\, \left(\frac{b+1}{2}\right)_n }
{\left(\frac{c}{2}\right)_n\, \left(\frac{c+1}{2}\right)_n\, (1)_n}
\end{eqnarray*}
Using the formula (\ref{inteq6}) and using the fact that $\Gamma(a+1)= a\Gamma(a)$, the aforementioned equation reduces to
\begin{flushleft}
$\displaystyle \sum_{n=0}^{\infty} \frac{(n+1)^2 \,(a)_n\, \left(\frac{b}{2}\right)_n\, \left(\frac{b+1}{2}\right)_n }
{\left(\frac{c}{2}\right)_n\, \left(\frac{c+1}{2}\right)_n\, (1)_n}$
\end{flushleft}
\begin{eqnarray*}
&=& \left(\frac{(a)_2\, (b)_4}{(c-a-b-2)_{2}\,(c-a)_2}\right)\left(\frac{\Gamma(c)\, \Gamma(c-a-b)}{\Gamma(c-b)\,\Gamma(c-a)}\right)\\ 
&& \qquad\qquad \times \,  _{2}F_1(a+2,b+4;c-a+2;-1)\\ 
&& + 3\left(\frac{(a)\, (b)_2}{(c-a-b-1)\,(c-a)}\right)\left(\frac{\Gamma(c)\, \Gamma(c-a-b)}{\Gamma(c-b)\,\Gamma(c-a)}\right)\\ 
&& \qquad\qquad \times \, _{2}F_1(a+1,b+2;c-a+1;-1)\\ 
&& + \left(\frac{\Gamma(c)\, \Gamma(c-a-b)}{\Gamma(c-a)\, \Gamma(c-b)}\right)\, _{2}F_1(a,b;c-a;-1).
\end{eqnarray*}
Hence,\\
$\displaystyle\sum_{n=0}^{\infty} \frac{(n+1)^2(a)_n\, \left(\frac{b}{2}\right)_n\, \left(\frac{b+1}{2}\right)_n }
{\left(\frac{c}{2}\right)_n\, \left(\frac{c+1}{2}\right)_n\, (1)_n}$
\begin{eqnarray*}
&=&  \frac{\Gamma(c)\, \Gamma(c-a-b)}{\Gamma(c-a)\,\Gamma(c-b)} \,\bigg[\frac{(a)_2\,(b)_4}{(c-a)_2\,(c-a-b-2)_{2}\,}\,\cr 
&& \qquad\qquad \times \, _{2}F_1(a+2,b+4;c-a+2;-1)\,\cr 
&& +\, 3\left(\frac{a\,(b)_2}{(c-a)(c-a-b-1)}\right)\, _{2}F_1(a+1,b+2;c-a+1;-1)\cr 
&&\,+\, _{2}F_1(a,b;c-a;-1)\bigg],
\end{eqnarray*}
Which completes the proof of (2).\\

(3)   Using $(n+1)^3=n(n-1)(n-2)+6n(n-1)+7n+1$, we can write

\begin{flushleft}
$\displaystyle \sum_{n=0}^{\infty} \frac{(n+1)^3 \,(a)_n\, \left(\frac{b}{2}\right)_n\, \left(\frac{b+1}{2}\right)_n }
{\left(\frac{c}{2}\right)_n\, \left(\frac{c+1}{2}\right)_n\, (1)_n}$
\end{flushleft}
\begin{eqnarray*}
&=&\left( \frac{(a)_3\,(b)_6}{(c)_6}\right)\sum_{n=0}^{\infty} \frac{(a+3)_{n}\, \left(\frac{b}{2}+3\right)_{n}\, \left(\frac{b+1}{2}+3\right)_{n} }
{\left(\frac{c}{2}+3\right)_{n}\, \left(\frac{c+1}{2}+3\right)_{n}\, (1)_{n}}\\
&& \,+6\, \left(\frac{(a)_2\,(b)_4}{(c)_4}\right) \sum_{n=0}^{\infty} \frac{(a+2)_{n}\, \left(\frac{b}{2}+2\right)_{n}\, \left(\frac{b+1}{2}+2\right)_{n} }
{\left(\frac{c}{2}+2\right)_{n}\, \left(\frac{c+1}{2}+2\right)_{n}\, (1)_{n}}\\
&& \,+7\,\left( \frac{a\,(b)_2}{(c)_2}\right)\sum_{n=0}^{\infty} \frac{(a+1)_{n}\, \left(\frac{b}{2}+1\right)_{n}\, \left(\frac{b+1}{2}+1\right)_{n} }
{\left(\frac{c}{2}+1\right)_{n}\, \left(\frac{c+1}{2}+1\right)_{n}\, (1)_{n}} +\sum_{n=0}^{\infty} \frac{(a)_n\, \left(\frac{b}{2}\right)_n\, \left(\frac{b+1}{2}\right)_n }
{\left(\frac{c}{2}\right)_n\, \left(\frac{c+1}{2}\right)_n\, (1)_n}
\end{eqnarray*}
Using the formula (\ref{inteq6}) and using the fact that $\Gamma(a+1)= a\Gamma(a)$, the aforementioned equation reduces to
\begin{flushleft}
$\displaystyle\sum_{n=0}^{\infty} \frac{(n+1)^3(a)_n\, \left(\frac{b}{2}\right)_n\, \left(\frac{b+1}{2}\right)_n }
{\left(\frac{c}{2}\right)_n\, \left(\frac{c+1}{2}\right)_n\, (1)_n}$
\end{flushleft}
\begin{eqnarray*}
 &=&    \frac{\Gamma(c)\, \Gamma(c-a-b)}{\Gamma(c-a)\,\Gamma(c-b)} \,\bigg[\left(\frac{(a)_3\,(b)_6}{(c-a)_3\,(c-a-b-3)_{3}\,}\right)\cr
&&\qquad\qquad \times\, _{2}F_1(a+3,b+6;c-a+3;-1)\,\cr
&&+ 6\, \left(\frac{(a)_2\,(b)_4}{(c-a)_2\,(c-a-b-2)_{2}\,}\right)\,\, _{2}F_1(a+2,b+4;c-a+2;-1)\,\cr
&&+\, 7\, \left(\frac{a\,(b)_2}{(c-a)(c-a-b-1)}\right)\, _{2}F_1(a+1,b+2;c-a+1;-1)\cr
&&\,+\, _{2}F_1(a,b;c-a;-1)\bigg].
\end{eqnarray*}
which completes the proof.\\

(4)  Let $a\neq 1$, $b\neq 1,\, 2$ and $c > \max\{a+1, a+b-1\}$. It is found that
\begin{flushleft}
$\displaystyle\sum_{n=0}^{\infty} \frac{(a)_n\, \left(\frac{b}{2}\right)_n\, \left(\frac{b+1}{2}\right)_n }
{\, \left(\frac{c}{2}\right)_n\, \left(\frac{c+1}{2}\right)_n\, (1)_{n+1}}$
\end{flushleft}
\begin{eqnarray*}
&=& \left(\frac{\left(c-1\right)\, \left(c-2\right)}{(a-1)\, \left(b-1\right)\, \left(b-2\right)}\right)  \left[\sum_{n=0}^{\infty} \frac{(a-1)_n\, \left(\frac{b}{2}-1\right)_n\, \left(\frac{b+1}{2}-1\right)_n }
{\left(\frac{c}{2}-1\right)_n\, \left(\frac{c+1}{2}-1\right)_n\, (1)_{n}} -1\right] \\ \\
&=& \left(\frac{(c-a-1)\,(c-a-b)\,\Gamma(c)\Gamma(c-a-b)}{(a-1)\, (b-1)\, (b-2)\Gamma(c-b)\, \Gamma(c-a)}\right)\cr
&&  \qquad \qquad  \times \, _{2}F_1(a-1,b-2;c-a-1;-1) - \left(\frac{(c-2)\, (c-1)}{(a-1)\, (b-1)\, (b-2)}\right).
\end{eqnarray*}
Hence the desired result follows.
\epf
\end{sloppypar}
\section{Starlikeness of $z\, _3F_2\left(a,\,\frac{b}{2},\, \frac{b+1}{2};\,\frac{c}{2},\, \frac{c+1}{2};z\right)$}
\bthm\label{ch3thm1eq0}
 Let $a,\, b \in {\Bbb C} \backslash \{ 0 \} $,\, $c > 0$\,  and\, $c > |a|+|b|+1.$ A sufficient condition for the function $z\, _3F_2\left(a,\,\frac{b}{2},\, \frac{b+1}{2};\,\frac{c}{2},\, \frac{c+1}{2};z\right) $ to belong to the class $ \es^{*}_{\lambda}, \,  0 < \lambda  \leq 1 $ is that
\beq\label{ch3thm1eq1}
   \frac{\Gamma(c)\, \Gamma(c-|a|-|b|)}{\Gamma(c-|a|)\, \Gamma(c-|b|)}\left[\left(\frac{|a|\, (|b|)_2}{(c-|a|)\, (c -|a|-|b|-1)}\right)\right. \qquad\qquad \qquad \qquad \qquad\nonumber\\ \nonumber \\
\left.   \times\, _{2}F_1(|a|+1,|b|+2;c-|a|+1;-1) + \lambda\, _{2}F_1(|a|,|b|;c-|a|;-1)\right] \leq 2\lambda.
\eeq
\ethm
\bpf  Let $f(z)=z\, _3F_2\left(a,\,\frac{b}{2},\, \frac{b+1}{2};\,\frac{c}{2},\, \frac{c+1}{2};z\right)$, then, by the equation (\ref{inteq2}), it is enough to show that
\begin{eqnarray*}
  T &=& \sum_{n=2}^{\infty}(n+\lambda-1)|A_n|\leq \lambda.
\end{eqnarray*}
Using the fact $|(a)_n|\leq  (|a|)_n$, one can get
\begin{eqnarray*}
  T &\leq&  \sum_{n=0}^{\infty} \left((n+1)\,\frac{(|a|)_{n}\left(\frac{|b|}{2}\right)_{n}\, \left(\frac{|b|+1}{2}\right)_{n}}{\left(\frac{c}{2}\right)_{n}\, \left(\frac{c+1}{2}\right)_{n}(1)_{n}}\right)\\ \\
  && \qquad+(\lambda-1)\sum_{n=0}^{\infty} \left(\frac{(|a|)_{n}\left(\frac{|b|}{2}\right)_{n}\, \left(\frac{|b|+1}{2}\right)_{n}}{\left(\frac{c}{2}\right)_{n}\, \left(\frac{c+1}{2}\right)_{n}(1)_{n}}\right)-\lambda.
\end{eqnarray*}
Using (\ref{inteq6}) and the result (1) of Lemma \ref{ch3lem1eq1} in the aforesaid equation, we get
\begin{eqnarray*}
T &\leq& \frac{\Gamma(c)\, \Gamma(c-|a|-|b|)}{\Gamma(c-|a|)\, \Gamma(c-|b|)}\bigg[\left(\frac{|a|\, (|b|)_2}{(c-|a|)\, (c -|a|-|b|-1)}\right)\nonumber \qquad \qquad \cr \cr
&& \times\, _{2}F_1(|a|+1,|b|+2;c-|a|+1;-1)\, +\, \lambda\, _{2}F_1(|a|,|b|;c-|a|;-1)\bigg]-\lambda.
\end{eqnarray*}
Because of (\ref{ch3thm1eq1}), the above expression is bounded above by $\lambda$, and hence,
\begin{eqnarray*}
  T &\leq&  \frac{\Gamma(c)\, \Gamma(c-|a|-|b|)}{\Gamma(c-|a|)\, \Gamma(c-|b|)}\bigg[\left(\frac{|a|\, (|b|)_2}{(c-|a|)\, (c -|a|-|b|-1)}\right)\nonumber \qquad \qquad \cr \cr
&& \times\, _{2}F_1(|a|+1,|b|+2;c-|a|+1;-1)\, +\, \lambda\, _{2}F_1(|a|,|b|;c-|a|;-1)\bigg]-\lambda.
\end{eqnarray*}
Therefore, $z\, _3F_2\left(a,\,\frac{b}{2},\, \frac{b+1}{2};\,\frac{c}{2},\, \frac{c+1}{2};z\right) $ belongs to the class $\es^{*}_{\lambda}. $
\epf
\bthm\label{ch3thm2eq001}
 Let $a,\, b \in {\Bbb C} \backslash \{ 0 \},\, c > 0, \, |a|\neq1,\, |b| \neq 1,\, 2,$ and $c > \max\{|a|+1, |a|+|b|-1\}.$ For  $ 0 < \lambda \leq 1$ and $ 0 \leq \beta < 1$, assume that
 \beq\label{ch3thm2eq1}
   \frac{\Gamma(c)\,\Gamma(c-|a|-|b|)}{\Gamma(c-|a|)\, \Gamma(c-|b|)}\left[\left(\frac{(\lambda-1)\,(c-|a|-|b|)\, (c-|a|-1)}{(|a|-1)\,(|b|-1)\,(|b|-2)}\right)\, \qquad\qquad\qquad\qquad\qquad\right.\nonumber\\ \nonumber \\
   \left.\times _{2}F_1(|a|-1,|b|-2;c-|a|-1;-1)\,+ \, _{2}F_1(|a|,|b|;c-|a|;-1)\right] \qquad\nonumber \\ \nonumber\\
  \leq \lambda\left(1+\frac{1}{2(1-\beta)}\right)+\frac{(\lambda-1)\,(c-1)\,(c-2)}{(|a|-1)(|b|-1)(|b|-2)}.
\eeq
 Then, the integral  operator $\mathcal{I}^{a,\,\frac{b}{2},\, \frac{b+1}{2}}_{\frac{c}{2},\, \frac{c+1}{2}}(f)$ maps $ \mathcal{R}(\beta)$ into $\es^{*}_{\lambda}$.
\ethm
\bpf
Let $a,\, b \in {\Bbb C} \backslash \{ 0 \},\, c > 0, \, |a|\neq1,\, |b| \neq 1,\, 2,$ and $c > |a|+|b|-1.$ For  $ 0 < \lambda \leq 1$ and $ 0 \leq \beta < 1$.\\

Consider the integral operator $\mathcal{I}^{a,\,\frac{b}{2},\, \frac{b+1}{2}}_{\frac{c}{2},\, \frac{c+1}{2}}(f)(z)$ defined by (\ref{inteq7}). According to (\ref{inteq2}), we need to show that
\begin{eqnarray}\label{thm2eq002}
  T &=& \sum_{n=2}^{\infty}(n+\lambda-1)|A_n|\leq \lambda,
\end{eqnarray}
where $A_n$ is given by (\ref{inteq007}). Then, we have
\begin{eqnarray*}
  T &=&  \sum_{n=2}^{\infty} [n+(\lambda-1)]\, \left|\frac{(a)_{n-1}\left(\frac{b}{2}\right)_{n-1}\, \left(\frac{b+1}{2}\right)_{n-1}}{\left(\frac{c}{2}\right)_{n-1}\, \left(\frac{c+1}{2}\right)_{n-1}(1)_{n-1}}\right| |a_n|
\end{eqnarray*}
Using (\ref{inteq3})  in the aforementioned equation, we have
\begin{eqnarray*}
  T   &\leq&  2(1-\beta)\bigg[\sum_{n=0}^{\infty}\, \left(\frac{(|a|)_{n}\left(\frac{|b|}{2}\right)_{n}\, \left(\frac{|b|+1}{2}\right)_{n}}{\left(\frac{c}{2}\right)_{n}\, \left(\frac{c+1}{2}\right)_{n}(1)_{n}}\right)-1\cr \cr \cr
  &&\qquad + (\lambda-1)\sum_{n=0}^{\infty} \left(\frac{(|a|)_{n}\left(\frac{|b|}{2}\right)_{n}\, \left(\frac{|b|+1}{2}\right)_{n}}{\left(\frac{c}{2}\right)_{n}\, \left(\frac{c+1}{2}\right)_{n}(1)_{n}}\right) \left(\frac{1}{n+1}\right)  -(\lambda-1)\bigg]:=T_1.
\end{eqnarray*}
Using the formula (\ref{inteq6}) and the results (1) and (4) of Lemma \ref{ch3lem1eq1}, we find that
\begin{eqnarray*}
T_1&\leq& 2(1-\beta)\bigg[ \left(\frac{\Gamma(c)\,\Gamma(c-|a|-|b|)}{\Gamma(c-|a|)\, \Gamma(c-|b|)}\right)\bigg(\,\left(\frac{(\lambda-1)\,(c-|a|-|b|)\, \, (c-|a|-1)}{(|a|-1)\, (|b|-1)\, (|b|-2)}\right)\cr \cr
  &&  \times\, _{2}F_1(|a|-1,|b|-2;c-|a|-1;-1) \,+ \, _{2}F_1(|a|,|b|;c-|a|;-1) \bigg) \cr \cr
  &&  -\frac{(\lambda-1)\,(c-1)\,(c-2)}{(|a|-1)(|b|-1)(|b|-2)} -\lambda  \bigg].
\end{eqnarray*}
Under the condition (\ref{ch3thm2eq1})
\begin{eqnarray*}
2(1-\beta)\qquad\qquad\qquad\qquad\qquad\qquad\qquad\qquad\qquad\qquad\qquad\qquad\qquad\qquad&&\\
  \times\bigg[ \left(\frac{\Gamma(c)\,\Gamma(c-|a|-|b|)}{\Gamma(c-|a|)\, \Gamma(c-|b|)}\right)\bigg( \left(\frac{(\lambda-1)\,(c-|a|-|b|)\, \, (c-|a|-1)}{(|a|-1)\, (|b|-1)\, (|b|-2)}\,\right) &&\cr \cr
\times \,_{2}F_1(|a|-1,|b|-2;c-|a|-1;-1) \,+ \, _{2}F_1(|a|,|b|;c-|a|;-1) \bigg) && \cr \cr
 -\,\frac{(\lambda-1)\,(c-1)\,(c-2)}{(|a|-1)(|b|-1)(|b|-2)} -\lambda \bigg] &\leq& \lambda.
\end{eqnarray*}
Thus, we have the inequalities $T \leq T_1 \leq \lambda $,  and hence  (\ref{thm2eq002}) hold. Therefore, it is concluded that the operator $\mathcal{I}^{a,\,\frac{b}{2},\, \frac{b+1}{2}}_{\frac{c}{2},\, \frac{c+1}{2}}(f)$ maps $ \mathcal{R}(\beta)$ into $\es^{*}_{\lambda}$, which completes the proof of the theorem.
\epf

When, $\lambda =1 $, we get the following result from Theorem \ref{ch3thm2eq001}.
\bcor
Let $a,\, b \in {\Bbb C} \backslash \{ 0 \},\, c > 0,\, $ and $c > |a|+|b|.$ For  $ 0 \leq \beta < 1$. Assume that
 \beq\label{cor2eq1}
   \left(\frac{\Gamma(c)\,\Gamma(c-|a|-|b|)}{\Gamma(c-|a|)\, \Gamma(c-|b|)}\right)\, _{2}F_1(|a|,|b|;c-|a|;-1) \leq 1+\frac{1}{2(1-\beta)}.\nonumber
\eeq
 Then, the integral  operator $\mathcal{I}^{a,\,\frac{b}{2},\, \frac{b+1}{2}}_{\frac{c}{2},\, \frac{c+1}{2}}(f)$ maps $ \mathcal{R}(\beta)$ into $\es^{*}_{1}$.
\ecor

\bthm\label{ch3thm3eq0}  Let $a,\, b \in {\Bbb C} \backslash \{ 0 \} $,\, $c > 0$\, and $c > |a|+|b|+2.$ For  $ 0 < \lambda \leq 1$. If
 \beq\label{ch3thm3eq1}
\left(\frac{\Gamma(c)\,\Gamma(c-|a|-|b|)}{\Gamma(c-|a|)\, \Gamma(c-|b|)}\right)\, \bigg[ \left(\frac{(|a|)_2\,(|b|)_4}{(c-|a|)_2\, (c-|a|-|b|-2)_{2}}\,\right)\qquad\qquad\qquad\qquad \cr \cr
\times \,_{2}F_1(|a|+2,|b|+4;c-|a|+2;-1)\qquad\qquad\,
\cr \cr \,+ \, \left(\frac{(\lambda+2)\, (|a|)\,(|b|)_2}{(c-|a|)\, (c-|a|-|b|-1)}\right)\, _{2}F_1(|a|+1,|b|+2;c-|a|+1;-1)\qquad
  \cr \cr +\qquad\, \lambda \, _{2}F_1(|a|,|b|;c-|a|;-1)\bigg] \leq 2\lambda,
\eeq
then the integral  operator $\mathcal{I}^{a,\,\frac{b}{2},\, \frac{b+1}{2}}_{\frac{c}{2},\, \frac{c+1}{2}}(f)$ maps $\mathcal{S}$ to $\es^{*}_{\lambda}$.
\ethm
\bpf Let $a,\, b \in {\Bbb C} \backslash \{ 0 \} $,\, $c > 0$\, $c > |a|+|b|+2$  and  $ 0 < \lambda \leq 1$. \\

Suppose that the integral  operator $\mathcal{I}^{a,\,\frac{b}{2},\, \frac{b+1}{2}}_{\frac{c}{2},\, \frac{c+1}{2}}(f)(z)$ is defined by (\ref{inteq7}).
In view of (\ref{inteq2}), it is enough to show that
\begin{eqnarray*}
  T &=& \sum_{n=2}^{\infty}(n+\lambda-1)|A_n|\leq \lambda.
\end{eqnarray*}
where $A_n$ is given by (\ref{inteq007}).
Using the fact $|(a)_n|\leq  (|a|)_n$ and the equation ($\ref{inteq00}$) in the aforementioned equation, it is derived that
\begin{eqnarray*}
  T &\leq & \sum_{n=2}^{\infty} n\,(n+(\lambda-1)) \left(\frac{(|a|)_{n-1}\left(\frac{|b|}{2}\right)_{n-1}\, \left(\frac{|b|+1}{2}\right)_{n-1}}{\left(\frac{c}{2}\right)_{n-1}\, \left(\frac{c+1}{2}\right)_{n-1}(1)_{n-1}}\right)
\end{eqnarray*}
Using (1) and (2) of Lemma \ref{ch3lem1eq1}, it is find that
\begin{eqnarray*}
T &\leq& \frac{\Gamma(c)\,\Gamma(c-|a|-|b|)}{\Gamma(c-|a|)\, \Gamma(c-|b|)}\, \bigg[ \left(\frac{(|a|)_2\,(|b|)_4}{(c-|a|)_2\, (c-|a|-|b|-2)_{2}}\,\right)\cr \cr
&&\qquad\qquad\times \,_{2}F_1(|a|+2,|b|+4;c-|a|+2;-1)\,\cr \cr
&&\,+ \, \left(\frac{(\lambda+2)\, |a| \,(|b|)_2}{(c-|a|)\, (c-|a|-|b|-1)}\right)\, _{2}F_1(|a|+1,|b|+2;c-|a|+1;-1) \cr \cr
&&\, +\, \lambda\, _{2}F_1(|a|,|b|;c-|a|;-1)\bigg] -\lambda.
\end{eqnarray*}
By (\ref{ch3thm3eq1}), the above expression is bounded above by $\lambda$, and hence,
\begin{eqnarray*}
 \frac{\Gamma(c)\,\Gamma(c-|a|-|b|)}{\Gamma(c-|a|)\, \Gamma(c-|b|)}\, \bigg[ \left(\frac{(|a|)_2\,(|b|)_4}{(c-|a|)_2\, (c-|a|-|b|-2)_{2}}\,\right)\qquad\qquad\qquad \qquad\qquad&&\cr \cr
\qquad\qquad\times \,_{2}F_1(|a|+2,|b|+4;c-|a|+2;-1)\,
 &&\cr \cr \,+ \, \left(\frac{(\lambda+2)\, (|a|)\,(|b|)_2}{(c-|a|)\, (c-|a|-|b|-1)}\right)\, _{2}F_1(|a|+1,|b|+2;c-|a|+1;-1)
 && \cr \cr \qquad +\, \lambda \, _{2}F_1(|a|,|b|;c-|a|;-1)\bigg] -\lambda &\leq& \lambda.
\end{eqnarray*}
Under the stated condition, the integral operator $\mathcal{I}^{a,\,\frac{b}{2},\, \frac{b+1}{2}}_{\frac{c}{2},\, \frac{c+1}{2}}(f)(z)$ maps $\es$ into $\es^{*}_{\lambda}$.
\epf
\section{Convexity of $z\, _3F_2\left(a,\,\frac{b}{2},\, \frac{b+1}{2};\,\frac{c}{2},\, \frac{c+1}{2};z\right)$ }
\bthm\label{ch3thm10eq1}
 Let $a,\, b \in {\Bbb C} \backslash \{ 0 \} $,\, $c > 0$,\, $c > |a|+|b|+2$  and  $ 0 < \lambda \leq 1$. A sufficient condition for the function $z\, _3F_2\left(a,\,\frac{b}{2},\, \frac{b+1}{2};\,\frac{c}{2},\, \frac{c+1}{2};z\right)$ to belong to the class $ \mathcal{C}_{\lambda}$ is that
\beq\label{ch3thm10eq10}
\left(\frac{\Gamma(c)\,\Gamma(c-|a|-|b|)}{\Gamma(c-|a|)\, \Gamma(c-|b|)}\right)\, \bigg[ \left(\frac{(|a|)_2\,(|b|)_4}{(c-|a|)_2\, (c-|a|-|b|-2)_{2}}\,\right) \qquad\qquad\qquad\qquad\qquad\cr \cr
\qquad\times \,_{2}F_1(|a|+2,|b|+4;c-|a|+2;-1)\,\qquad\qquad
 \cr \cr \,+ \, \left(\frac{(\lambda+2)\, (|a|)\,(|b|)_2}{(c-|a|)\, (c-|a|-|b|-1)}\right)\, _{2}F_1(|a|+1,|b|+2;c-|a|+1;-1)
  \qquad\qquad \cr \cr \qquad\qquad+\, \lambda\, _{2}F_1(|a|,|b|;c-|a|;-1)\bigg] \leq 2\lambda.\nonumber
\eeq
\ethm
\bpf The proof is similar to Theorem \ref{ch3thm3eq0}. So we omit the details.
\epf
\bthm\label{ch3thm11eq0}  Let $a,\, b \in {\Bbb C} \backslash \{ 0 \} $,\, $c > 0$,\, $c > |a|+|b|+1$  and  $ 0 < \lambda \leq 1$.  For $0 \leq \beta <1 $, it is assumed that
 \beq\label{ch3thm11eq1}
   \left(\frac{\Gamma(c)\,\Gamma(c-|a|-|b|)}{\Gamma(c-|a|)\, \Gamma(c-|b|)}\right)  \qquad \qquad \qquad \qquad \qquad \qquad \qquad \qquad \qquad \qquad \qquad \qquad \qquad \nonumber\\ \times \bigg[ \left(\frac{(|a|)\,(|b|)_2}{(c-|a|)\, (c-|a|-|b|-1)}\right) \, _{2}F_1(|a|+1,|b|+2;c-|a|+1;-1) \qquad \qquad \cr \,+\, \lambda \, _{2}F_1(|a|,|b|;c-|a|;-1)\bigg]  \leq \lambda\left( \frac{1}{2(1-\beta)}+1\right).\nonumber
\eeq
Then, the operator $\mathcal{I}^{a,\,\frac{b}{2},\, \frac{b+1}{2}}_{\frac{c}{2},\, \frac{c+1}{2}}(f)$ maps $\mathcal{R}(\beta)$ into $ \mathcal{C}_{\lambda} $.
\ethm
\bpf The proof is similar to Theorem \ref{ch3thm2eq001}. So we omit the details.
\epf
\bthm\label{thm12eq0}  Let $a,\, b \in {\Bbb C} \backslash \{ 0 \} $,\, $c > 0$,\, $c > |a|+|b|+3$  and  $ 0 < \lambda \leq 1$. If
 \beq\label{ch3thm12eq1}
 &&\left(\frac{\Gamma(c)\,\Gamma(c-|a|-|b|)}{\Gamma(c-|a|)\, \Gamma(c-|b|)}\right)\,\nonumber\cr \cr
 && \qquad\times \bigg[ \left(\frac{(|a|)_3\,(|b|)_6}{(c-|a|)_3\, (c-|a|-|b|-3)_{3}}\,\right)\,\, _{2}F_1(|a|+3,|b|+6;c-|a|+3;-1)\,\cr \cr
 && \qquad +\, \left(\frac{(\lambda+5)\,(|a|)_2\,(|b|)_4}{(c-|a|)_2\, (c-|a|-|b|-2)_{2}}\,\right)\, \,_{2}F_1(|a|+2,|b|+4;c-|a|+2;-1)\,\cr \cr
 && \qquad+ \, \left(\frac{(3\lambda+4)\, (|a|)\,(|b|)_2}{(c-|a|)\, (c-|a|-|b|-1)}\right)\,\, _{2}F_1(|a|+1,|b|+2;c-|a|+1;-1)\cr \cr
 && \qquad\qquad\qquad\qquad\qquad\qquad\qquad\qquad\qquad\qquad\qquad+\, \lambda\, _{2}F_1(|a|,|b|;c-|a|;-1)\bigg] \leq 2\lambda,\nonumber
\eeq
then, $\mathcal{I}^{a,\,\frac{b}{2},\, \frac{b+1}{2}}_{\frac{c}{2},\, \frac{c+1}{2}}(f)$ maps $\es$ into $ \mathcal{C}_{\lambda} $.
\ethm
\bpf Let $a,\, b \in {\Bbb C} \backslash \{ 0 \} $,\, $c > 0$,\, $c > |a|+|b|+3$  and  $ 0 < \lambda \leq 1$.\\

Suppose the integral operator $\mathcal{I}^{a,\,\frac{b}{2},\, \frac{b+1}{2}}_{\frac{c}{2},\, \frac{c+1}{2}}(f)(z)$ is defined by (\ref{inteq7}). In view of the sufficient condition given in (\ref{inteq}), it is enough to prove that
\begin{eqnarray*}
  T &=& \sum_{n=2}^{\infty}\,n\, (n+\lambda-1)\, |A_n|\leq \lambda .
\end{eqnarray*}
i.e.,
\begin{eqnarray*}
  T &=& \sum_{n=2}^{\infty}n\, (n+\lambda-1)\, \left|\left(\frac{(a)_{n-1}\left(\frac{b}{2}\right)_{n-1}\, \left(\frac{b+1}{2}\right)_{n-1}}{\left(\frac{c}{2}\right)_{n-1}\, \left(\frac{c+1}{2}\right)_{n-1}(1)_{n-1}}\right)\right|\, |a_n|\leq \lambda.
\end{eqnarray*}
Using the fact that $|(a)_n|\leq  (|a|)_n$ and $(\ref{inteq00})$ in the aforementioned equation, it is derived
\begin{eqnarray*}
  T &\leq&\sum_{n=0}^{\infty} \frac{(n+1)^3\,(|a|)_{n-1}\left(\frac{|b|}{2}\right)_{n}\, \left(\frac{|b|+1}{2}\right)_{n}}{\left(\frac{c}{2}\right)_{n}\, \left(\frac{c+1}{2}\right)_{n}(1)_{n}}\\ \\
  &&+(\lambda-1)\sum_{n=0}^{\infty} \frac{(n+1)^2\,(|a|)_{n}\left(\frac{|b|}{2}\right)_{n}\, \left(\frac{|b|+1}{2}\right)_{n}}{\left(\frac{c}{2}\right)_{n}\, \left(\frac{c+1}{2}\right)_{n}(1)_{n}}-\lambda.
\end{eqnarray*}
Using (2) and (3) of Lemma \ref{ch3lem1eq1}, we find that
\begin{eqnarray*}
T &\leq & \left(\frac{\Gamma(c)\,\Gamma(c-|a|-|b|)}{\Gamma(c-|a|)\, \Gamma(c-|b|)}\right)\\
 && \times\, \bigg[ \left(\frac{(|a|)_3\,(|b|)_6}{(c-|a|)_3\, (c-|a|-|b|-3)_{3}}\,\right) \,_{2}F_1(|a|+3,|b|+6;c-|a|+3;-1)\,\cr \cr
 && \,+\left(\frac{(\lambda+5)\,(|a|)_2\,(|b|)_4}{(c-|a|)_2\, (c-|a|-|b|-2)_{2}}\,\right) \,_{2}F_1(|a|+2,|b|+4;c-|a|+2;-1)\,\cr \cr
 && \,+ \, \left(\frac{\,(3\lambda+4)\, (|a|)\,(|b|)_2}{(c-|a|)\, (c-|a|-|b|-1)}\right)\, _{2}F_1(|a|+1,|b|+2;c-|a|+1;-1) \cr \cr
 &&+\, \lambda\, _{2}F_1(|a|,|b|;c-|a|;-1)\bigg] -\lambda.
\end{eqnarray*}
By the equation (\ref{ch3thm12eq1}), the above expression is bounded above by $\lambda$, and hence,
\begin{eqnarray*}
\left(\frac{\Gamma(c)\,\Gamma(c-|a|-|b|)}{\Gamma(c-|a|)\, \Gamma(c-|b|)}\right)\qquad\qquad\qquad\qquad\qquad\qquad\qquad\qquad\qquad\qquad\qquad&&\\
 \times\, \bigg[ \left(\frac{(|a|)_3\,(|b|)_6}{(c-|a|)_3\, (c-|a|-|b|-3)_{3}}\,\right) \,_{2}F_1(|a|+3,|b|+6;c-|a|+3;-1)&&\,\cr \cr
 \,+\left(\frac{(\lambda+5)\,(|a|)_2\,(|b|)_4}{(c-|a|)_2\, (c-|a|-|b|-2)_{2}}\,\right) \,_{2}F_1(|a|+2,|b|+4;c-|a|+2;-1)&&\,\cr \cr
 \,+ \, \left(\frac{\,(3\lambda+4)\, (|a|)\,(|b|)_2}{(c-|a|)\, (c-|a|-|b|-1)}\right)\, _{2}F_1(|a|+1,|b|+2;c-|a|+1;-1) &&\cr 
 +\, \lambda\, _{2}F_1(|a|,|b|;c-|a|;-1)\bigg] -\lambda &\leq& \lambda.
\end{eqnarray*}
Hence, the integral operator $\mathcal{I}^{a,\,\frac{b}{2},\, \frac{b+1}{2}}_{\frac{c}{2},\, \frac{c+1}{2}}(f)(z)$ maps $\es$ into $\mathcal{C}_{\lambda}$ and the proof is complete.
\epf
\section{Admissibility condition of $z\,  _3F_2\left(a,\,\frac{b}{2},\, \frac{b+1}{2};\,\frac{c}{2},\, \frac{c+1}{2};z\right)$ in $UCV$.}
\bthm\label{ch3thm7eq1}
 Let $a,\, b \in {\Bbb C} \backslash \{ 0 \} $,\, $c > 0$\, and $c > |a|+|b|+2$.  A sufficient condition for the function $z\, _3F_2\left(a,\,\frac{b}{2},\, \frac{b+1}{2};\,\frac{c}{2},\, \frac{c+1}{2};z\right) $ to belong to the class ${\UCV}$ is that\\
 \beq\label{ch3thm7eq10}
   \left(\frac{\Gamma(c)\,\Gamma(c-|a|-|b|)}{\Gamma(c-|a|)\, \Gamma(c-|b|)}\right)\qquad\qquad \qquad\qquad\qquad\qquad\qquad\qquad\qquad\qquad\qquad\nonumber&&\\ \times\bigg[ \, \left(\frac{2 \,(|a|)_2\,(|b|)_4}{(c-|a|)_2\, (c-|a|-|b|-2)_{2}}\,\right) \,_{2}F_1(|a|+2,|b|+4;c-|a|+2;-1)\,&& \cr \cr \,+\, 5\, \left(\frac{(|a|)\,(|b|)_2}{(c-|a|)\, (c-|a|-|b|-1)}\right)\, _{2}F_1(|a|+1,|b|+2;c-|a|+1;-1) &&  \cr \cr +\, _{2}F_1(|a|,|b|;c-|a|;-1)\bigg] &\leq& 2.
   \eeq
\ethm
\bpf
Let $a,\, b \in {\Bbb C} \backslash \{ 0 \} $,\, $c > 0$\, and  $c > |a|+|b|+2$.\\

Let  $\displaystyle f(z)=z\,  _3F_2\left(a,\,\frac{b}{2},\, \frac{b+1}{2};\,\frac{c}{2},\, \frac{c+1}{2};z\right)$. Then, by (\ref{lem4eq1}), it is enough to show that
\begin{eqnarray*}
  T &=& \sum_{n=2}^{\infty}\, n\, (2n-1)\,|A_n|\leq 1.
\end{eqnarray*}
where $A_n$ is given by (\ref{inteq007}). Using the fact $|(a)_n|\leq  (|a|)_n$,
\begin{eqnarray*}
  T &\leq& 2\sum_{n=0}^{\infty}  \left(\frac{(n+1)^2\,(|a|)_{n}\left(\frac{|b|}{2}\right)_{n}\, \left(\frac{|b|+1}{2}\right)_{n}}{\left(\frac{c}{2}\right)_{n}\, \left(\frac{c+1}{2}\right)_{n}(1)_{n}}\right)
  -\sum_{n=0}^{\infty} \left(\frac{ (n+1)\, (|a|)_{n}(|b|)_{n}(c)_{n}}{(|b|+1)_{n}(c+1)_{n}(1)_{n}}\right)-1.
\end{eqnarray*}
Using (1) and (2) of Lemma \ref{ch3lem1eq1} in the aforementioned  equation, we find that
\begin{eqnarray*}
  T &\leq& \frac{\Gamma(c)\,\Gamma(c-|a|-|b|)}{\Gamma(c-|a|)\, \Gamma(c-|b|)}\, \bigg[ \, \left(\frac{2 \,(|a|)_2\,(|b|)_4}{(c-|a|)_2\, (c-|a|-|b|-2)_{2}}\,\right) \\ 
  &&\qquad\qquad\times\,_{2}F_1(|a|+2,|b|+4;c-|a|+2;-1)\,\cr 
 && \,+\, 5\, \left(\frac{ (|a|)\,(|b|)_2}{(c-|a|)\, (c-|a|-|b|-1)}\right)\, _{2}F_1(|a|+1,|b|+2;c-|a|+1;-1) \cr 
 && +\, _{2}F_1(|a|,|b|;c-|a|;-1)\bigg]-1.
\end{eqnarray*}
Because of (\ref{ch3thm7eq10}), the above expression is bounded above by 1, and hence,
\begin{eqnarray*}
\left(\frac{\Gamma(c)\,\Gamma(c-|a|-|b|)}{\Gamma(c-|a|)\, \Gamma(c-|b|)}\right)\qquad\qquad \qquad\qquad\qquad\qquad\qquad\qquad\qquad\qquad\qquad\nonumber&&\\ \times\bigg[ \, \left(\frac{2 \,(|a|)_2\,(|b|)_4}{(c-|a|)_2\, (c-|a|-|b|-2)_{2}}\,\right) \,_{2}F_1(|a|+2,|b|+4;c-|a|+2;-1)\,&& \cr \cr \,+\, 5\, \left(\frac{(|a|)\,(|b|)_2}{(c-|a|)\, (c-|a|-|b|-1)}\right)\, _{2}F_1(|a|+1,|b|+2;c-|a|+1;-1) && \cr \qquad\qquad +\, _{2}F_1(|a|,|b|;c-|a|;-1)\bigg]-1 &\leq& 1.
\end{eqnarray*}
Therefore, $z\, _3F_2\left(a,\,\frac{b}{2},\, \frac{b+1}{2};\,\frac{c}{2},\, \frac{c+1}{2};z\right) $ belongs to the class $ UCV. $
\epf
\bthm\label{ch3thm8eq0}  Let $a,\, b \in {\Bbb C} \backslash \{ 0 \} $,\, $c > 0$,\, $c > |a|+|b|+1$ and $0 \leq \beta <1 $. Assume that
 \beq\label{ch3thm8eq1}
\frac{\Gamma(c)\,\Gamma(c-|a|-|b|)}{\Gamma(c-|a|)\, \Gamma(c-|b|)}\,\, \bigg[ \left(\frac{2\,(|a|)\,(|b|)_2}{(c-|a|)\, (c-|a|-|b|-1)}\right) \, \qquad\qquad\qquad\cr \cr
 \times_{2}F_1(|a|+1,|b|+2;c-|a|+1;-1)\, \, +\, _{2}F_1(|a|,|b|;c-|a|;-1)\bigg] \nonumber\\  \leq \frac{1}{2(1-\beta)}+1.
\eeq
Then, $\mathcal{I}^{a,\,\frac{b}{2},\, \frac{b+1}{2}}_{\frac{c}{2},\, \frac{c+1}{2}}(f)$  maps $\mathcal{R}(\beta)$ into ${\UCV}$.
\ethm
\bpf Let $a,\, b \in {\Bbb C} \backslash \{ 0 \} $,\, $c > 0$,\, $c > |a|+|b|+1$  and  $ 0 < \beta \leq 1$.

Then consider the integral operator $\mathcal{I}^{a,\,\frac{b}{2},\, \frac{b+1}{2}}_{\frac{c}{2},\, \frac{c+1}{2}}(f)$ given in (\ref{inteq7}). According to sufficient condition given in (\ref{lem4eq1}), it is enough to show that
\begin{eqnarray*}
  T &:=& \sum_{n=2}^{\infty}n\, (2n-1)\,|A_n|\leq 1,
\end{eqnarray*}
where $A_n$ is given by (\ref{inteq007}). Using the fact $|(a)_n|\leq  (|a|)_n$ and $(\ref{inteq3})$ in the aforementioned equation, it is found that
\begin{eqnarray*}
  T &\leq& 2(1-\beta)\sum_{n=2}^{\infty} n\,(2n-1)\, \left(\frac{(|a|)_{n-1}\left(\frac{|b|}{2}\right)_{n-1}\, \left(\frac{|b|+1}{2}\right)_{n-1}}{\left(\frac{c}{2}\right)_{n-1}\, \left(\frac{c+1}{2}\right)_{n-1}(1)_{n-1}\, n}\right)
\end{eqnarray*}
Using the formula (\ref{inteq6}) and (1) of Lemma \ref{ch3lem1eq1}, it is derived that
\begin{eqnarray*}
T &\leq& 2(1-\beta)\bigg[ \frac{\Gamma(c)\,\Gamma(c-|a|-|b|)}{\Gamma(c-|a|)\, \Gamma(c-|b|)}\, \bigg[ \left(\frac{2\,(|a|)\,(|b|)_2}{(c-|a|)\, (c-|a|-|b|-1)}\right)\,\cr \cr  &&\, \times\, _{2}F_1(|a|+1,|b|+2;c-|a|+1;-1)\,+\, _{2}F_1(|a|,|b|;c-|a|;-1)\bigg] -1 \bigg].
\end{eqnarray*}
By (\ref{ch3thm8eq1}), the aforementioned expression is bounded above by $1$, and hence,
\begin{eqnarray*}
&& 2(1-\beta)\bigg[ \frac{\Gamma(c)\,\Gamma(c-|a|-|b|)}{\Gamma(c-|a|)\, \Gamma(c-|b|)}\, \bigg[ \left(\frac{2\,(|a|)\,(|b|)_2}{(c-|a|)\, (c-|a|-|b|-1)}\right)\,\cr \cr  &&\qquad\, \times\, _{2}F_1(|a|+1,|b|+2;c-|a|+1;-1)\,+\, _{2}F_1(|a|,|b|;c-|a|;-1)\bigg] -1 \bigg] \leq 1.
\end{eqnarray*}

Therefore, the operator $\mathcal{I}^{a,\,\frac{b}{2},\, \frac{b+1}{2}}_{\frac{c}{2},\, \frac{c+1}{2}}(f)(z)$ maps $\mathcal{R}(\beta)$ into $UCV$, and the result follows.
\epf
\section{Inclusion Properties of $z\,  _3F_2\left(a,\,\frac{b}{2},\, \frac{b+1}{2};\,\frac{c}{2},\, \frac{c+1}{2};z\right) $  in $\es_p$-CLASS}
\bthm\label{ch3thm4eq0}
 Let $a,\, b \in {\Bbb C} \backslash \{ 0 \} $,\, $c > 0$\, and $c > |a|+|b|+1$.  A sufficient condition for the function $z\, _3F_2\left(a,\,\frac{b}{2},\, \frac{b+1}{2};\,\frac{c}{2},\, \frac{c+1}{2};z\right) $ to belong to the class $\es_p$ is that
 \beq\label{ch3thm4eq1}
   \left(\frac{\Gamma(c)\,\Gamma(c-|a|-|b|)}{\Gamma(c-|a|)\, \Gamma(c-|b|)}\right)\,\,  \bigg[ \left(\frac{2\,(|a|)\,(|b|)_2}{(c-|a|)\, (c-|a|-|b|-1)}\right)\qquad\qquad\qquad\qquad\qquad\, \cr \cr \, \times\, _{2}F_1(|a|+1,|b|+2;c-|a|+1;-1)\, +\, _{2}F_1(|a|,|b|;c-|a|;-1)\bigg] \leq 2.\nonumber
\eeq
\ethm
\bpf
The proof is similar to Theorem \ref{ch3thm8eq0}. So we omit the details.
\epf
\bthm\label{ch3thm5eq0}  Let $a,\, b \in {\Bbb C} \backslash \{ 0 \},  c > 0,\,   |a| \neq 1 ,\, |b|\neq 1,\, 2$,\, $c > \max\{|a|+1, |a|+|b|-1\}$ and $0 \leq \beta < 1$. Assume that
\beq\label{ch3thm5eq1}
&&\frac{\Gamma(c)\,\Gamma(c-|a|-|b|)}{\Gamma(c-|a|)\, \Gamma(c-|b|)}\,  \bigg[ 2\, _{2}F_1(|a|,|b|;c-|a|;-1)\,\cr \cr
 &&\qquad +\left(\frac{(c-|a|-1)\,(c-|a|-|b|)}{(|a|-1)(|b|-1)(|b|-2)}\right)\, _{2}F_1(|a|-1,|b|-2;c-|a|-1;-1)\,\bigg]\,\,
  \cr \cr && \qquad\qquad\qquad\qquad\qquad\qquad\qquad \qquad\,+ \frac{(c-1)\,(c-2)}{(|a|-1)\,(|b|-1)\,(|b|-2)} \leq \frac{1}{2(1-\beta)}+1.
\eeq
Then, $\mathcal{I}^{a,\,\frac{b}{2},\, \frac{b+1}{2}}_{\frac{c}{2},\, \frac{c+1}{2}}(f)$ maps $\mathcal{R}(\beta)$ into $\es_p$ class.
\ethm
\bpf Let $a,\, b \in {\Bbb C} \backslash \{ 0 \} ,  c > 0,\,   |a| \neq 1 ,\, |b|\neq 1,\, 2$,\, $c > \max\{|a|+1, |a|+|b|-1\}$ and $0 \leq \beta < 1$. \\
Consider the integral operator $\mathcal{I}^{a,\,\frac{b}{2},\, \frac{b+1}{2}}_{\frac{c}{2},\, \frac{c+1}{2}}(f)$ given by (\ref{inteq7}). In the view of (\ref{lem2eq1}), it is enough to show that
\begin{eqnarray*}
  T &:=& \sum_{n=2}^{\infty}(2n-1)|A_n|\leq 1,
\end{eqnarray*}
where $A_n$ is given by (\ref{inteq007}). It is proven that
\begin{eqnarray*}
  T &\leq & \sum_{n=2}^{\infty}(2n-1)\left(\frac{(|a|)_{n-1}\left(\frac{|b|}{2}\right)_{n-1}\, \left(\frac{|b|+1}{2}\right)_{n-1}}{\left(\frac{c}{2}\right)_{n-1}\, \left(\frac{c+1}{2}\right)_{n-1}(1)_{n-1}}\right)\, |a_n|\leq 1.
\end{eqnarray*}
Using the inequality $|(a)_n|\leq  (|a|)_n$ and $(\ref{inteq3})$ in the aforesaid equation, it is derived that
\begin{eqnarray*}
T  &\leq & 2(1-\beta)\bigg[2\sum_{n=0}^{\infty}  \frac{(n+1)(|a|)_{n}\left(\frac{|b|}{2}\right)_{n}\, \left(\frac{|b|+1}{2}\right)_{n}}{\left(\frac{c}{2}\right)_{n}\, \left(\frac{c+1}{2}\right)_{n}(1)_{n+1}}\cr \cr
&&-\sum_{n=0}^{\infty}  \frac{(|a|)_{n}\left(\frac{|b|}{2}\right)_{n}\, \left(\frac{|b|+1}{2}\right)_{n}}{\left(\frac{c}{2}\right)_{n}\, \left(\frac{c+1}{2}\right)_{n}(1)_{n+1}}-1\bigg].
\end{eqnarray*}
Using the equation (\ref{inteq6}) and (4) of Lemma \ref{ch3lem1eq1}, it is found that
\begin{eqnarray*}
  T &\leq&  2(1-\beta)\bigg[ \frac{\Gamma(c)\,\Gamma(c-|a|-|b|)}{\Gamma(c-|a|)\, \Gamma(c-|b|)}\, \bigg[ 2\, _{2}F_1(|a|,|b|;c-|a|;-1)\, \cr \cr
  &&\  -\left(\frac{(c-|a|-1)\,(c-|a|-|b|)}{(|a|-1)(|b|-1)(|b|-2)}\right)\, \, _{2}F_1(|a|-1,|b|-2;c-|a|-1;-1)\,\bigg]\cr \cr
  && \,+\, \frac{(c-1)\,(c-2)}{(|a|-1)\,(|b|-1)\,(|b|-2)}-1  \bigg].
\end{eqnarray*}
By the condition (\ref{ch3thm5eq1}), the aforementioned expression is bounded above by 1, and hence,
\begin{eqnarray*}
  2(1-\beta)\bigg[ \frac{\Gamma(c)\,\Gamma(c-|a|-|b|)}{\Gamma(c-|a|)\, \Gamma(c-|b|)}\, \bigg[ 2\, \, _{2}F_1(|a|,|b|;c-|a|;-1)
  \qquad\qquad\qquad\qquad&&\cr \cr \,-\left(\frac{(c-|a|-1)\,(c-|a|-|b|)}{(|a|-1)(|b|-1)(|b|-2)}\right)\, _{2}F_1(|a|-1,|b|-2;c-|a|-1;-1)\,\bigg]\,&&\cr \cr
  +\, \frac{(c-1)\,(c-2)}{(|a|-1)\,(|b|-1)\,(|b|-2)}-1  \bigg] &\leq& 1.
\end{eqnarray*}
Under the stated condition, the operator $\mathcal{I}^{a,\,\frac{b}{2},\, \frac{b+1}{2}}_{\frac{c}{2},\, \frac{c+1}{2}}(f)(z)$ maps $\mathcal{R}(\beta)$ into $\es_p$ and  the proof is complete.
\epf
\bthm\label{ch3thm6eq0} Let $a,\, b \in {\Bbb C} \backslash \{ 0 \} $,\, $c > 0$\, and $c > |a|+|b|+2$. Suppose $a,\, b$, and $c$ satisfy the condition
 \beq\label{ch3thm6eq1}
\left(\frac{\Gamma(c)\,\Gamma(c-|a|-|b|)}{\Gamma(c-|a|)\, \Gamma(c-|b|)}\right)\qquad\qquad\qquad\qquad\qquad\qquad\qquad\qquad\qquad\qquad\qquad\qquad&& \cr \cr
\qquad\qquad\times \bigg[ \, \left(\frac{2 \,(|a|)_2\,(|b|)_4}{(c-|a|)_2\, (c-|a|-|b|-2)_{2}}\,\right) \,_{2}F_1(|a|+2,|b|+4;c-|a|+2;-1)\,&&\cr \cr
\,+\, 5\, \left(\frac{(|a|)\,(|b|)_2}{(c-|a|)\, (c-|a|-|b|-1)}\right)\, _{2}F_1(|a|+1,|b|+2;c-|a|+1;-1)&& \cr \cr
+\, _{2}F_1(|a|,|b|;c-|a|;-1)\bigg] &\leq& 2.
\eeq
Then, $\mathcal{I}^{a,\,\frac{b}{2},\, \frac{b+1}{2}}_{\frac{c}{2},\, \frac{c+1}{2}}(f)$ maps $\es$ into $ \es_{p} $ class.
\ethm
\bpf
The proof is similar to Theorem \ref{ch3thm7eq1}. So we omit the details.
\epf


\begin{thebibliography}{150}
\bibitem{Anbu-Parva-2000}
        {M. Anbu Durai and R. Parvatham,} {\it Convolutions with Hypergeometric Functions,} Bull. Malaysian Math. Sc. Soc. (Second Series) {\bf 23} (2000) 153-161.

\bibitem{Chandru-prabha-2019}
        {K.Chandrasekran\ and\ D.J.Prabhakaran}
        {\it Geometric Properties of Clausen's  Hypergeometric Functions,}
        Preprint.

\bibitem{Driver-Johnston-2006}
        {K. A. Driver\ and\ S. J. Johnston,}
        {\it  An integral representation of some hypergeometric functions,} Electron. Trans. Numer. Anal. {\bf 25} (2006), 115--120.

\bibitem{Peter-L-Duren-book-1983}
         {P. L. Duren}, {\it Univalent functions}
         (Grundlehren der   mathematischen Wissenschaften 259, New York, Bcerlin, Heidelberg,Tokyo), Springer--Verlag, (1983).

\bibitem{A-W-Goodman-1983-book}
         {A. W. Goodman}, {\it Univalent functions}, Vol.I and Vol.II, Tampa Florida Mariner Publishing Company, (1983).

\bibitem{Good-1991-Ann-PM}
        {A. W. Goodman},
        {\it On uniformly convex functions},
         Ann. Polon. Math. {\bf 56} (1991), no.~1, 87--92.

\bibitem{Good-1991-JMAA}
        {A. W. Goodman},
        {\it On uniformly starlike functions},
        J. Math. Anal. Appl. {\bf 155} (1991), no.~2, 364--370.

\bibitem{Ma-Minda-1992-Ann-PM}
        { W. C. Ma\ and\ D. Minda},
        {\it Uniformly convex functions},
         Ann. Polon. Math. {\bf 57} (1992), no.~2, 165--175.

\bibitem{MacGregor-1962-Trans-ams}
        {T. H. MacGregor},
        {\it Functions whose derivative has a positive real part},
        Trans. Amer. Math. Soc. {\bf 104} (1962), 532--537.

\bibitem{Parva-prabha-2001-Far-East}
        {R. Parvatham} and {D.J. Prabhakaran},
        {\it On the Hohlov convolution operator of the class $S_p,$}
         Far East J. Math.Sci.Special volume, (2001), 217--228.



\bibitem{Ponnu-saba-1997}
        {S. Ponnusamy\ and\ S. Sabapathy,}
        {Geometric properties of generalized hypergeometric functions,}
         Ramanujan J. {\bf 1} (1997), no.~2, 187--210.

\bibitem{Robertson-1936}
        {Robertson, Malcolm I. S.}
        {\it On the theory of univalent functions}.
         Ann. of Math. (2) 37 (1936), no. 2, 374--408.

\bibitem{Ronn-1993-Proc-ams}
        { F. R\o nning},
        {\it Uniformly convex functions and a corresponding class of starlike functions},
         Proc. Amer. Math. Soc. {\bf 118} (1993), no.~1, 189--196.

\bibitem{Subram-Murugu-1995}
        {K. G. Subramanian\ et al.,}
        {\it Subclasses of uniformly convex and uniformly starlike functions,}
         Math. Japon. {\bf 42} (1995), no.~3, 517--522.

\bibitem{Subra-Sudharsan-1998}
        {K. G. Subramanian\ et al.,}
        {\it Classes of uniformly starlike functions,}
        Pub. Math. Debrecen {\bf 53} (1998), no.~3-4, 309--315.
\end{thebibliography}
\end{document}